\PassOptionsToPackage{unicode}{hyperref}
\PassOptionsToPackage{hyphens}{url}
\documentclass[
]{article}
\usepackage{amsmath,amssymb}
\usepackage{lmodern}
\usepackage{iftex}
\ifPDFTeX
  \usepackage[T1]{fontenc}
  \usepackage[utf8]{inputenc}
  \usepackage{textcomp} 
\else 
  \usepackage{unicode-math}
  \defaultfontfeatures{Scale=MatchLowercase}
  \defaultfontfeatures[\rmfamily]{Ligatures=TeX,Scale=1}
\fi
\IfFileExists{upquote.sty}{\usepackage{upquote}}{}
\IfFileExists{microtype.sty}{
  \usepackage[]{microtype}
  \UseMicrotypeSet[protrusion]{basicmath} 
}{}
\makeatletter
\@ifundefined{KOMAClassName}{
  \IfFileExists{parskip.sty}{%
    \usepackage{parskip}
  }{
    \setlength{\parindent}{0pt}
    \setlength{\parskip}{6pt plus 2pt minus 1pt}}
}{
  \KOMAoptions{parskip=half}}
\makeatother
\usepackage{xcolor}
\usepackage{graphicx}
\graphicspath{ {C:/Users/andre/Dropbox/Probability/Random graphs/Parking on trees with power law tails - my article/Slika1.png} }
\makeatletter
\def\maxwidth{\ifdim\Gin@nat@width>\linewidth\linewidth\else\Gin@nat@width\fi}
\def\maxheight{\ifdim\Gin@nat@height>\textheight\textheight\else\Gin@nat@height\fi}
\makeatother
\setkeys{Gin}{width=\maxwidth,height=\maxheight,keepaspectratio}
\makeatletter
\def\fps@figure{htbp}
\makeatother
\ifLuaTeX
  \usepackage{selnolig}  
\fi

\def\title#1{{\Large\bf  \begin{center} #1 \vspace{0pt} \end{center}  } }
\def\authors#1{{\large\bf \begin{center} #1 \vspace{0pt} \end{center} } }
\def\university#1{{\sl \begin{center} #1 \vspace{0pt} \end{center} } }
\def\inst#1{\unskip$^{#1}$}

\IfFileExists{bookmark.sty}{\usepackage{bookmark}}{\usepackage{hyperref}}
\IfFileExists{xurl.sty}{\usepackage{xurl}}{} 
\hypersetup{
  hidelinks,
  pdfcreator={LaTeX via pandoc}}

\author{}
\date{}

\begin{document}

\title{Short survey of results and open problems for parking problems on random trees}

\bigskip

%
%

\authors{Andrej Srakar\inst{1}
     }

\smallskip

%
%

 \university{
 \inst{1} Institute for Economic Research and University of Ljubljana, Slovenia, andrej.srakar@ier.si}

\bigskip

\noindent {\large\bf Abstract}

\medskip

Parking problems derive from works in combinatorics by Konheim and Weiss
in the 1960s. In a memorable contribution, Lackner and Panholzer (2016)
studied parking on a random tree and established a phase transition for
this process when \(m \approx \frac{n}{2}\). This relates to the
renowned result by David Aldous of convergence results on Erdős-Renyi
random graphs of order \(n^{\frac{2}{3}}\). In a series of recent
articles, Contat and coauthors have studied the problem in various
random tree contexts and derived several novel scaling limit and phase
transition results. We survey the present state-of-the-art of this
literature and point to its extensions, open directions and
possibilities, in particular related to the study of problem in
different metric topologies. My intent it to point to importance of this
line of research and novel open problems for future study.



\begin{enumerate}
\def\labelenumi{\arabic{enumi}.}
\item
  Introduction
\end{enumerate}

Parking functions were introduced by Konheim and Weiss (1966) in their
investigations of a linear probing collision resolution scheme for hash
tables. Since then, they have attracted plenty of attention and proven
to be a fertile source of interesting mathematics. They have found
extensions to research in representation theory (Pak and Postnikov,
1994; Armstrong et al., 2015), polytopes (Stanley and Pitman, 2002), the
sandpile model (Cori and Le Borgne, 2003), probability theory and
stochastic processes (Lackner and Panholzer, 2016; Contat and Curien,
2023), and the theory of Macdonald polynomials in combinatorics (Haiman,
1994), as just few examples.

They have been translated to a probabilistic problem by Lackner and
Panholzer in their article from 2016 in Journal of Combinatorial Theory,
Series A. In a series of recent articles they received many novel
probabilistic results. Intent of my short survey contribution is to
resume the present state-of-the-art on addressing the topic and point to
several open problems and interesting future extensions.

\begin{enumerate}
\def\labelenumi{\arabic{enumi}.}
\setcounter{enumi}{1}
\item
  Parking functions in combinatorics
\end{enumerate}

Kovalinka and Towari explain basic features of parking functions in
combinatorics, and point to extensions to a subclass of rational parking
functions. As explained in their contribution (Kovalinka and Towari,
2021), an integer sequence \((x_{1},\ldots,x_{n})\) is a parking
function if its weakly increasing rearrangement \((z_{1},\ldots,z_{n})\)
satisfies \(0 \leq z_{i} \leq i - 1\) for \(i = 1,\ldots n\). This
definition implies that rearranging the entries in one parking function
results in another. Haiman (1994) was the first to study the \(S_{n}\)
action on the set of parking functions of length \(n\). Two decades
later, Berget and Rhoades (2014) studied the following seemingly
unrelated representation \(\sigma_{n}\) of \(S_{n}\). Let \(K_{n}\)
denote the complete graph with vertex set
\(\lbrack n\rbrack = \{ 1,\ldots,n\}\). Given a subgraph
\(G \subseteq K_{n}\), we attach to it the polynomial
\(p(G) = \prod_{ij \in E(G)}^{}(x_{i} - x_{j}\mathbb{) \in C\lbrack}x_{1},\ldots,x_{n}\rbrack\).Here
\(E(G)\) refers to the set of edges of \(G\) and we record those by
listing the smaller number first. Define \(V_{n}\) to be the
\(\mathbb{C}\)-linear span of \(p(G)\) over all \(G\) for which the
complement \(\overline{G}\) is a connected graph. \(V_{n}\) first
appears in the work of Postnikov and Shapiro (2004). The natural action
of \(S_{n}\) on \(\mathbb{C\lbrack}x_{1},\ldots,x_{n}\rbrack\) that
permutes variables gives an action on \(V_{n}\) because relabeling
vertices preserves connectedness. Berget and Rhoades (2014, Theorem 2)
also established the remarkable fact that the restriction of
\(\sigma_{n}\) to \(S_{n - 1}\) is isomorphic to \(\rho_{n - 1}\).

For \(n \geq 1\), we denote by \(\mathbb{Z}_{n}\) the set of integers
modulo \(n\). Typically, representatives from residue classes modulo
\(n\) will be implicitly assumed to belong to \(\{ 0,\ldots,n - 1\}\).
In the below, let \(S_{n}\) denote the symmetric group consisting of
permutations of \(\lbrack n\rbrack\). We use both the cycle notation and
one-line notation for permutations. Within the latter let \(\pi_{i}\)
denote the image of \(i\) under the permutation \(\pi\) for a positive
integer \(i\).

A partition \(\lambda = (\lambda_{1},\ldots,\lambda_{l})\) is a weakly
decreasing sequence of positive integers. The
\(\lambda_{i}\)\textquotesingle s\textquotesingle{} are the parts of
\(\lambda\), their sum its size, and their number its length, which is
denoted by \(l(\lambda)\). If \(\lambda\) has size \(n\), then we denote
this by \(\lambda \vdash n\). Furthermore, letting \(m_{i}\) denote the
multiplicity of the part \(i\) in \(\lambda\) for \(i \geq 1\), we set
\(z_{\lambda} = \prod_{i \geq 1}^{}i^{m_{i}}m_{i}!\). The cycle type of
a permutation \(\pi\) is a partition that we denote \(\lambda(\pi)\).

We consider the following distinguished bases for the ring of symmetric
functions \(\Lambda\): the power sum symmetric functions
\(\{ p_{\lambda}:\lambda \vdash n\}\), the complete homogeneous
symmetric functions \(\{ h_{\lambda}:\lambda \vdash n\}\), and the Schur
symmetric functions \(\{ s_{\lambda}:\lambda \vdash n\}\).

Representation theory of the symmetric group is intimately tied to
\(\Lambda\) and the connection is made explicit by the Frobenius
characteristic. Given a representation \(\rho\) of \(S_{n}\), denote the
corresponding character by \(\chi_{\rho}\). Then

\[Frob(\rho) = \frac{1}{n!}\sum_{\pi \in S_{n}}^{}{\chi_{\rho}(\pi)}p_{\lambda(\pi)} = \sum_{\lambda \vdash n}^{}{\chi_{\rho}(\lambda)}\frac{p_{\lambda}}{z_{\lambda}}\]

Under \(Frob\), the irreducible representation of \(S_{n}\)
corresponding to the partition \(\mu \vdash n\) gets mapped to the Schur
function \(s_{\mu}\). As a special case, we have the equality
\(\sum_{\lambda \vdash n}^{}z_{\lambda}^{- 1}p_{\lambda} = h_{n}\).

An integer sequence \((x_{1},\ldots,x_{n})\) is a parking function if
its weakly increasing rearrangement \((z_{1},\ldots,z_{n})\) satisfies
\(0 \leq z_{i} \leq i - 1\) for \(i = 1,\ldots,n\). We denote by
\({PF}_{n}\) the set of all parking functions of length \(n\). For
example,

\[{PF}_{2} = \left\{ 00,01,10 \right\},\]

\[{PF}_{3} = \left\{ 000,001,010,100,002,020,200,011,101,110,012,021,102,120,201,210 \right\},\]

and the weakly increasing elements of \({PF}_{4}\) are
\(0000,0001,0011,0111,0002,0012,0112,0022,0122,0003,0013,0113,0023,\)
and \(0123\). Observe that there are 14 such elements in \({PF}_{4}\).
More generally, we have that the number of weakly increasing elements in
\({PF}_{n}\) is the \(n\)th Catalan number
\({Cat}_{n} = \frac{1}{n + 1}\begin{pmatrix}
2n \\
n \\
\end{pmatrix}\). It is well known that
\(\left| {PF}_{n} \right| = {(n + 1)}^{n - 1}\). This is seen through
the following result in Foata and Riordan (1974, where it is attributed
to H. O. Pollak):

Theorem (Pollak, in Foata in Riorda, 1974)). The map
\({PF}_{n} \rightarrow \mathbb{Z}_{n + 1}^{n - 1}\), given by

\[(x_{1},\ldots,x_{n}) \mapsto (x_{2} - x_{1},\ldots,x_{n} - x_{n - 1})\]

where subtraction is performed modulo \(n + 1\), is a bijection.

For a partition \(\lambda = (\lambda_{1},\ldots,\lambda_{l}) \vdash n\)
the number of fixed points of the action of the permutation with cycle
decomposition
\((1,\ldots,\lambda_{1})(\lambda_{1} + 1,\ldots,\lambda_{1} + \lambda_{2})\cdots\)
is equal to the number of sequences
\((\alpha_{1},\ldots,\alpha_{n - 1}) \in \mathbb{Z}_{n + 1}^{n - 1}\)
satisfying \(\alpha_{i} = 0\) for
\(i \in \lbrack n - 1\rbrack \smallsetminus \{\lambda_{1},\lambda_{1} + \lambda_{2},\ldots,\lambda_{1} + \ldots + \lambda_{l - 1}\}\).
It follows that the character \(\chi_{\rho_{n}}\) of \(\rho_{n}\)
satisfies

\[\chi_{\rho_{n}}(\pi) = {(n + 1)}^{l - 1},\]

where \(l = l(\lambda(\pi))\).

In their original article, Konheim and Weiss consider the structure of
systems for filing, cataloguing and storing units of information, where
each record book or information unit has a natural name or record
identification number associated with it. The set of all possible names
\(\{ a_{1},a_{2},\ldots,a_{m}\}\) is usually very large in comparison to
the actual number \(r\) of records
\(\{ a_{i_{1}},a_{i_{2}},\ldots,a_{i_{r}}\}\) that are to be stored in
any one problem. The storage procedure consists of assigning to each
record \(a_{i_{k}}\) a unique record location number
\(A_{i_{k}} \in \{ 0,1,\ldots,n - 1\}\) where \(n\) is the size of the
storage and \(r \leq n\). Typical values of \(m\) and \(n\) are
\(2^{36}\) and \(2^{10}\) respectively. The problem is to devise a
procedure for assigning the record location numbers so that the time
needed to store and recover a record, knowing only its name, is
minimized.

In most situations, \(\{ a_{i_{1}},a_{i_{2}},\ldots,a_{i_{r}}\}\) lacks
a definite structure and \(m\) is much larger than \(n\). Various
schemes for storage have been considered. One has been described by
Peterson (1957) as follows. One begins by randomly selecting a function
\(g:\{ a_{1},a_{2},\ldots,a_{m}\} \rightarrow \{ 0,1,\ldots,n - 1\}\).
The record location numbers
\(\{ A_{i_{1}},A_{i_{2}},\ldots,A_{i_{r}}\}\) of the records
\(a_{i_{1}},a_{i_{2}},\ldots,a_{i_{r}}\) are defined inductively as
follows:

\begin{enumerate}
\def\labelenumi{(\roman{enumi})}
\item
  \(A_{i_{1}} = g\left( a_{i_{1}} \right),\)
\item
  \(A_{i_{k}} = g\left( a_{i_{k}} \right) + s_{k}(modulo\ n),\)
\end{enumerate}

where \(s_{k}\) is the smallest nonnegative integer such that
\(g\left( a_{i_{k}} \right) + s_{k}(modulo\ n) \notin \{ A_{i_{1}},A_{i_{2}},\ldots,A_{i_{k - 1}}\}\).
To recover the record \(a_{i_{k}}\) one computes in succession the
record location numbers
\(g\left( a_{i_{k}} \right),g\left( a_{i_{k}} \right) + 1(modulo\ n),\ldots\),
comparing after each computation the name of the record stored in each
of the these locations with \(a_{i_{k}}\). Number of comparisons needed
to recover the record \(a_{i_{k}}\) is just \(s_{k} + 1\).

Konheim and Weiss establish some preliminaries. Let \(n\) be a positive
integer and

\[\pi = \left( \pi_{1},\pi_{2},\ldots,\pi_{n} \right) \in P_{n}\]

a permutation of the integers \(1,2,\ldots,n\). Define
\(\tau_{j,n},\tau_{n}\) and \(T_{n}\) by

\[\tau_{j,n}(\pi) = \max{\{ k:k \leq j,}\pi_{j} \geq \pi_{m}\ for\ m = j,j - 1,\ldots,j - k + 1\}\]

\[\tau_{n}(\pi) = \prod_{j = 1}^{n}{\tau_{j,n}(\pi)},\]

and

\[T_{n} = \sum_{\pi \in P_{n}}^{}{\tau_{n}(\pi)}\]

Then it holds that \(T_{n} = {(n + 1)}^{n - 1},n = 1,2,\ldots\) (Lemma
1, Konheim and Weiss, 1966).

Consider \(r\) balls \(B_{1},B_{2},\ldots,B_{r}\) which are to be placed
into \(n\) cells \(C_{0},C_{1},\ldots,C_{n - 1}\). We assume
\(r \leq n\). The location of the \(r\) balls are determined according
to the following occupancy discipline: suppose \(r\) fictitious cell
numbers \(\left( j_{1},j_{2},\ldots,j_{k},\ldots,j_{r} \right)\) have
been selected \((0 \leq j_{k} < n,1 \leq k \leq r)\). The actual
location of the \(k\)th ball \(B_{k}\), say \(l_{k}\), is defined
inductively according to the rules

\begin{enumerate}
\def\labelenumi{(\roman{enumi})}
\item
  \(l_{1} = j_{1}\),
\item
  for \(k \geq 2\), \(l_{k} = j_{k} + s_{k}(modulo\ n)\), where
  \(s_{k}\) is the smallest nonnegative integer such that
\end{enumerate}

\[l_{k} = j_{k} + s_{k}(modulo\ n) \notin \{ l_{1},l_{2},\ldots,l_{k - 1}\}\]

We let \(A\) denote the transformation

\[A:j = \left( j_{1},j_{2},\ldots j_{r} \right) \rightarrow Aj = 1 = \left( l_{1},\ldots l_{r} \right),\]

and set

\[\mathfrak{A}_{1} = \left\{ j:Aj = 1 \right\}.\]

Then, Konheim and Weiss prove the following lemma:

Lemma 2 (Konheim and Weiss, 1966): \(f(n,\ r) = n^{r - 1}(n - r)\)

Let
\(X = \left\{ a_{1},a_{2},\ldots,a_{m} \right\},Y = 0,1,2,\ldots,n - 1\}\)
and \(\mathcal{G}(X,Y) = \{ g:g:X \rightarrow Y\}\). The elements of
\(X\) are record identification numbers and the elements of \(Y\) are
record location numbers. Let
\(S = \{ a_{i_{1}},a_{i_{2}},\ldots,a_{i_{r}}\}\) be a fixed (ordered)
set of record identification numbers with \(r \leq n\). An element
\(g \in \mathcal{G}(X,Y)\) determines record location numbers for \(S\)
according to the rules:

\begin{enumerate}
\def\labelenumi{(\roman{enumi})}
\item
  the record location number for \(a_{i_{1}}\) is
  \(A_{i_{1}} = g(a_{i_{1}})\),
\item
  for \(k \geq 2\), the record location number for \(a_{i_{k}}\) is
  \(A_{i_{k}} = g\left( a_{i_{k}} \right) + s_{k}(modulo\ n)\) where
  \(s_{k}\) is the smallest nonnegative integer such that
  \(g\left( a_{i_{k}} \right) + s_{k}(modulo\ n) \notin \notin \{ A_{i_{1}},A_{i_{2}},\ldots,A_{i_{k - 1}}\}\).
\end{enumerate}

Let \((\Omega\mathcal{,E,}\mathbb{P})\) be a probability space. Let
\(G\) be a \(\mathcal{G}(X,Y)\)-valued random variable on
\((\Omega\mathcal{,E,}\mathbb{P})\) with
\(\mathbb{P}\left\{ \omega:G(\omega) = g \right\} = n^{- m},\ g \in \mathcal{G}(X,Y)\).

Konheim and Weiss prove the following theorems:

Theorem 1 (Konheim and Weiss, 1966):
\(\mathbb{P}\left\{ \omega:s_{k}(\omega) = j \right\} = \frac{1}{n^{k - 1}}\sum_{q = j}^{k - 1}\begin{pmatrix}
k - 1 \\
q \\
\end{pmatrix}{(q + 1)}^{q - 1} \bullet (n - k)(n - q - 1)^{k - q - 2},\)
\(\mathbb{E}\left\{ s_{k} \right\} = \frac{n - k}{2n^{k - 1}}\sum_{q = j}^{k - 1}\begin{pmatrix}
k - 1 \\
q \\
\end{pmatrix}{(q + 1)}^{q}q(n - q - 1)^{k - 2 - q}\).

Theorem 2 (Konheim and Weiss, 1966): Let \(\mu \in (0,1)\). Then

\[\lim_{n \rightarrow \infty}{\mathbb{E}\left\{ s_{\mu n} \right\}} = \frac{1}{2}\mu\frac{(2 - \mu)}{{(1 - \mu)}^{2}}.\]

In the final chapter of their article, Konheim and Weiss translate this
to parking problems. They define \(st.\) as a street with \(p\) parking
places. A car occupied by a man and his dozing wife enters \(st.\) at
the left and moves towards the right. The wife awakens at a capricious
moment and orders her husband to park immediately. He dutifully parks at
his present location, if it is empty, and if not, continues to the right
and parks at the next available space. If no space is available he
leaves \(st.\)

Suppose \(st.\) to be initially empty and \(c\) cars arrive with
independently capricious wives in each car. Konheim and Weiss calculate
the probability that they all find parking places. If by »capricious« is
meant that the probability of awakening in front of the \(i\)th parking
place is \(\frac{1}{p},\ 1 \leq i \leq p\), then the desired probability
is just

\[\mathbb{P}(c,p) = \frac{f(p + 1,c)}{p^{c}} = \left( 1 + \frac{1}{p} \right)^{c}\left( 1 - \frac{c}{p + 1} \right).\]

In particular it holds that
\(\lim_{p \rightarrow \infty}{\mathbb{P}(\mu,p,p)} = (1 - \mu)e^{\mu},\ 0 < \mu \leq 1.\)

The right hand side of the above expression for \(\mathbb{P}(c,p)\) is
also the probability that \(c\) cars will succeed in parking in \(st.\)
of length \(p\) (initially vacant) under the following more complicated
parking discipline: when the \(i\)the car stops he parks if the space is
free. If the space is occupied he performs a chance experiment; with
probability \(q_{i}\) he moves backward and with probability
\(1 - q_{i}\) he moves forward, in both cases seeking the first free
space.

\begin{enumerate}
\def\labelenumi{\arabic{enumi}.}
\setcounter{enumi}{2}
\item
  Parking on a random tree in Lackner and Panholzer
\end{enumerate}

In a 2016 article, Marie-Louise Lackner and Alois Panholzer (Lackner and
Panholzer, 2016) studied parking problems on a random tree. By this they
have put parking problem in a probability context and significantly
extended previous analyses of parking problems in combinatorics and
mathematics in general.

To explain their work, the following notation will prove useful. Given
an \(n\)-mapping \(f\), we define a binary relation \(\preccurlyeq_{f}\)
on \(\lbrack n\rbrack\) via

\[i \preccurlyeq_{f}j: \Longleftrightarrow \exists k\mathbb{\in N:}f^{k}(i) = j.\]

Thus \(i \preccurlyeq_{f}j\) holds if there exists a directed path from
\(i\) to \(j\) in the functional digraph \(G_{f}\), and we say that
\(j\) is a successor of \(i\) or that \(i\) is a predecessor of \(j\).
In this context a one-way street represents a total order, a tree
represents a certain partial order, where the root node is the maximal
element and a mapping represents a certain pre-order, i.e. a binary
relation that is transitive and reflexive.

The combinatorial structure of the functional digraph \(G_{f}\) of an
arbitrary mapping function \(f\) is well known: the weakly connected
components of \(G_{f}\) are cycles of rooted labelled trees. That is,
each connected component consists of rooted labelled trees whose root
nodes are connected by directed edges such that they form a cycle. We
call a node \(j\) for which there exists a \(k \geq 1\) such that
\(f^{k}(j) = j\), a cyclic node.

For ordinary parking functions it holds that changing the order of the
elements of a sequence does not affect its properts of being a parking
function or not. This can easily be generalized to parking functions for
mappings, and is resumed in the following lemma:

Lemma 2.1. (Lackner and Panholzer, 2016): A function
\(s:\lbrack m\rbrack \rightarrow \lbrack n\rbrack\) is a parking
function for a mapping
\(f:\lbrack n\rbrack \rightarrow \lbrack n\rbrack\) if and only if
\(s \circ \sigma\) is a parking function for \(f\) for any permutation
\(\sigma\) on \(\lbrack m\rbrack\).

Now let\textquotesingle s turn to parking functions where the number of
drivers does not coincide with the number of parking spaces. It is
well-known that a parking sequence
\(s:\lbrack m\rbrack \rightarrow \lbrack n\rbrack\) on a one-way street
is a parking function if and only if

\[\left| \left\{ k \in \lbrack m\rbrack:s_{k} \geq j \right\} \right| \leq n - j + 1,\ for\ all\ j \in \lbrack n\rbrack\]

In the above, the path \(s_{j} = y_{j} \rightsquigarrow \pi_{s}(j)\)
denotes the parking path of the \(j\)-th driver of \(s\) in the mapping
graph \(G_{f}\) starting with the preferred parking space \(s_{j}\) and
ending with the parking position \(\pi_{s}(j)\).

This can be generalized to \((n,m)\)-tree parking functions as follows.
It is known that a parking sequence
\(s:\lbrack m\rbrack \rightarrow \lbrack n\rbrack\) on a one-way street
is a parking function if and only if
\(\left| \left\{ k \in \lbrack m \right\rbrack:\ s_{k} \geq j \right\}| \leq n - j + 1,\ for\ all\ j \in \lbrack n\rbrack\).

Lemma 2.3. (Lackner and Panholzer, 2016): Given a rooted labelled tree
\(T\) of size \(|T| = n\) and a sequence \(s \in \lbrack n\rbrack^{m}\).
Then \(s\) is a tree parking function for \(T\) if and only if

\[\left| \left\{ k \in \lbrack m\rbrack:s_{k} \in T' \right\} \right| \leq \left| T^{'} \right|,\ for\ all\ subtrees\ T^{'}of\ T\ containing\ root(T)\]

Lackner and Panholzer estimate the number of parking functions. Given an
\(n\)-mapping \(f:\lbrack n\rbrack \rightarrow \lbrack n\rbrack\), let
us denote by \(S(f,m)\) the number of parking functions
\(s \in \lbrack n\rbrack^{m}\) for \(f\) with \(m\) drivers. Let \(T\)
be a rooted labelled tree. Tight bounds for \(S(T,m)\) are obtained by
them as follows:

Theorem 2.6. (Lackner and Panholzer, 2016): Let \({star}_{n}\) be the
rooted labelled tree of size \(n\) with root node \(n\) and the nodes
\(1,2,\ldots,\ n - 1\) attached to it. Furthermore, let \({chain}_{n}\)
be the rooted labelled tree of size \(n\) with root node \(n\) and node
\(j\) attached to node \((j + 1)\), for \(1 \leq j \leq n - 1\). Then,
for any rooted labelled tree \(T\) of size \(n\) it holds

\[S\left( {star}_{n},m \right) \leq S(T,m) \leq S\left( {chain}_{n},m \right),\]

yielding the bounds

\[n^{m} + \begin{pmatrix}
m \\
2 \\
\end{pmatrix}{(n - 1)}^{m - 1} \leq S(T,m) \leq (n - m + 1)(n + 1)^{m - 1},\ for\ 0 \leq m \leq n.\]

Lackner and Panholzer also study the total number of \((n,n)\)-mapping
parking functions \(M_{n} = M_{n,n}\), i.e. the number of pairs
\((f,s)\) with \(f \in M_{n}\) an \(n\)-mapping and
\(s \in \lbrack n\rbrack^{n}\) a parking sequence of length \(n\) for
the mapping \(f\), such that all drivers are successful. They derive the
following results.

Lemma 3.1. (Lackner and Panholzer, 2016): The total number \(C_{n}\) of
parking functions of length \(n\) for connected \(n\)-mappings is, for
\(n \geq 1\), given as follows:

\[C_{n} = n!(n - 1)!\sum_{j = 0}^{n - 1}\frac{{(2n)}^{j}}{j!}\]

Theorem 3.2. (Lackner and Panholzer, 2016): For all \(n \geq 1\) it
holds that the total numbers \(F_{n}\) and \(M_{n}\) of \((n,n)\)-tree
parking functions and \((n,n)\)-mapping parking functions, respectively,
satisfy:

\[M_{n} = n \bullet F_{n}.\]

Theorem 3.3. (Lackner and Panholzer, 2016): The total number \(M_{n}\)
of \((n,n)\)-mapping parking functions is for \(n \geq 1\) given as
follows:

\[M_{n} = n!(n - 1)!\sum_{j = 0}^{n - 1}\frac{{(n - j) \bullet (2n)}^{j}}{j!}.\]

Corollary 3.4. (Lackner and Panholzer, 2016): The total number \(F_{n}\)
of \((n,n)\)-tree parking functions is for \(n \geq 1\) given as
follows:

\[F_{n} = {((n - 1)!)}^{2}\sum_{j = 0}^{n - 1}\frac{{(n - j) \bullet (2n)}^{j}}{j!}.\]

Lackner and Panholzer also derive equivalent results and study the exact
and asymptotic behaviour of the total number of tree and mapping parking
functions for the general case of \(n\) parking spaces and
\(0 \leq m \leq n\) drivers. They analyze the total number \(F_{n,m}\)
of \((n,m)\)-tree parking functions, i.e. the number of pairs \((T,s)\),
with \(T \in \mathcal{T}_{n}\) a Cayley tree of size \(n\) and
\(s \in \lbrack n\rbrack^{m}\) a parking sequence of length \(m\) for
the tree \(T\), such that all drivers are successful. Furthermore,
\(F_{n,n} = F_{n}\) denotes the number of tree parking functions when
the number of parking spaces \(n\) coincides with the number of drivers
\(m\). They derive the following main results.

Theorem 4.4. (Lackner and Panholzer, 2016): For all \(n \geq 1\) it
holds that the total numbers \(F_{n,m}\) and \(M_{n,m}\) of
\((n,m)\)-tree parking functions and \((n,m)\)-mapping parking
functions, respectively, satisfy:

\[M_{n,m} = n \bullet F_{n,m}.\]

Theorem 4.5. (Lackner and Panholzer, 2016): The total number \(M_{n,m}\)
of \((n,m)\)-mapping parking functions is, for \(0 \leq m \leq n\) and
\(n \geq 1\), given as follows:

\[M_{n,m} = \frac{(n - 1)!m!n^{n - m}}{(n - m)!}\sum_{j = 0}^{m}{\begin{pmatrix}
2m - n - j \\
m - j \\
\end{pmatrix}\frac{{(n - j) \bullet (2n)}^{j}}{j!}}.\]

Corollary 4.6. (Lackner and Panholzer, 2016): The total number
\(F_{n,m}\) of \((n,m)\)-tree parking functions is, for
\(0 \leq m \leq n\) and \(n \geq 1\), given as follows:

\[F_{n,m} = \frac{(n - 1)!m!n^{n - m - 1}}{(n - m)!}\sum_{j = 0}^{m}{\begin{pmatrix}
2m - n - j \\
m - j \\
\end{pmatrix}\frac{{(n - j) \bullet (2n)}^{j}}{j!}}.\]

They derive two asymptotic results as follows.

Theorem 4.10. (Lackner and Panholzer, 2016): The total number
\(M_{n,m}\) of \((n,m)\)-mapping parking functions is asymptotically,
for \(n \rightarrow \infty\), given as folows (where \(\delta\) denotes
an arbitrary small, but fixed, constant):

\[M_{n,m}\sim\left\{ \begin{matrix}
\frac{n^{n + m + \frac{1}{2}}\sqrt{n - 2m}}{n - m}\ for\ 1 \leq m \leq \left( \frac{1}{2} - \delta \right)n \\
\frac{\sqrt{2}3^{\frac{1}{6}}\Gamma\left( \frac{2}{3} \right)n^{\frac{3n}{2} - \frac{1}{6}}}{\sqrt{\pi}}\ for\ m = \frac{n}{2} \\
\frac{m!}{(n - m)!} \bullet \frac{n^{2n - m + \frac{3}{2}}2^{2m - n + 1}}{(2m - n)^{\frac{5}{2}}}\ for\ \left( \frac{1}{2} + \delta \right)n \leq m \leq n \\
\end{matrix} \right.\ \]

Corollary 4.11. (Lackner and Panholzer, 2016): The probability
\(p_{n,m}\) that a randomly chosen pair \((f,s)\) with \(f\) an
\(n\)-mapping and \(s\) a sequence in \(\lbrack{n\rbrack}^{m}\),
represents a parking function is asymptotically, for
\(n \rightarrow \infty\) and \(m = \rho n\) with \(0 < \rho < 1\) fixed,
given as follows:

\[p_{n,m}\sim\left\{ \begin{matrix}
C_{<}(\rho)\ for\ 0 \leq \rho \leq \frac{1}{2} \\
C_{\frac{1}{2}} \bullet n^{- \frac{1}{6}}\ for\ \rho = \frac{1}{2} \\
C_{>}(\rho) \bullet n^{- 1} \bullet \left( D_{>}(\rho) \right)^{n}\ for\ \frac{1}{2} < \rho < 1 \\
\end{matrix} \right.\ \]

with

\[C_{<}(\rho) = \frac{\sqrt{1 - 2\rho}}{1 - \rho}\]

\[C_{\frac{1}{2}} = \sqrt{\frac{6}{\pi}}\frac{\Gamma(\frac{2}{3})}{3^{\frac{1}{3}}} \approx 1.298\ldots\]

\[C_{>}(\rho) = 2 \bullet \sqrt{\frac{\rho}{(1 - \rho){(2\rho - 1)}^{5}}}\]

\[D_{>}(\rho) = \left( \frac{4\rho}{e^{2}} \right)^{\rho}\frac{e}{2{(1 - \rho)}^{1 - \rho}}\]

Lackner and Panholzer also list the following open problems for research
in parking on random trees in the future:

\begin{enumerate}
\def\labelenumi{(\arabic{enumi})}
\item
  Given a tree \(T\) or a mapping \(f\), is it possible in general to
  give some simple characterization of the numbers \(S(T,m)\) and
  \(S(f,m)\), respectively?
\item
  With the approach presented, one can also study the total number of
  parking functions for other important tree families as, e.g., labelled
  binary trees or labelled ordered trees.
\item
  The problem of determining the total number of parking functions seems
  to be interesting for so-called increasing (or decreasing) tree
  families. For so-called recursive trees, i.e., unordered increasing
  trees, the approach presented could be applied, but the differential
  equations occurring do not seem to yield tractable solutions. For such
  tree families quantities such as the sums of parking functions as
  studied could be worthwhile treating as well.
\item
  As for ordinary parking functions one could analyse important
  quantities for tree and mapping parking functions. E.g., the so-called
  total displacement (which is of particular interest in problems
  related to hashing algorithms), i.e., the total driving distance of
  the drivers, or individual displacements (the driving distance of the
  k-th driver) seem to lead to interesting questions.
\item
  A refinement of parking functions can be obtained by studying what has
  been called defective parking functions or overflow, i.e., pairs
  \((T,s)\) or \((f,s)\), such that exactly \(k\) drivers are
  unsuccessful. Preliminary studies indicate that the approach presented
  is suitable to obtain results in this direction as well.
\item
  One could consider enumeration problems for some restricted parking
  functions for trees (or mappings).
\item
  Let us denote by \(X_{n}\) the random variable measuring the number of
  parking functions \(s\) with \(n\) drivers for a randomly chosen
  labelled unordered tree \(T\) of size \(n\). Then, due to our previous
  results, we get the expected value of \(X_{n}\) via
  \(\mathbb{E}\left( X_{n} \right) = \frac{F_{n}}{T_{n}}\sim\frac{\sqrt{2\pi}2^{n + 1}n^{n - \frac{1}{2}}}{e^{n}}\).
  However, with the approach presented here, it seems that we are not
  able to obtain higher moments or other results on the distribution of
  \(X_{n}\).
\end{enumerate}

\begin{enumerate}
\def\labelenumi{\arabic{enumi}.}
\setcounter{enumi}{3}
\item
  Recent interest of study of parking on random trees and its scaling
  limits
\end{enumerate}

Recently, Alice Contat has performed a lot of interesting work
addressing several open issues pointed already by Lackner and Panholzer
and proving many novel results. In her thesis work she dealed with the
study of parking models on random graphs and trees in a broad sense. She
investigated two algorithms to find a large independent set of a graph,
that is a subset of the vertices of the graph where no pair of vertices
are connected to each other. The first one uses a greedy procedure to
construct an independent set which is maximal for the inclusion order.
In the generic case, this subset has a positive density and she provided
example of large random graphs for which we can explicitly compute the
law of the size of the greedy independent set. Second is Karp--Sipser
algorithm which is optimal in the sense that there exists an independent
set with the maximal possible size which contains the subset of vertices
produced by algorithm.

She gave a precise localization of the phase transition for the
existence of a giant Karp--Sipser core for a configuration model with
vertices of degree 1, 2 and 3, and precisely analyzed its size at
criticality. Then, she examined the dynamical parking model introduced
by Konheim and Weiss on the line and considered a rooted tree where each
vertex represents a park spot and the edges are oriented towards the
root. She observed a phase transition and provided a localisation of it
for critical Bienaymé--Galton--Watson trees using the local limit, and
for the infinite binary trees via a combinatorial decomposition. On
critical trees, she also showed that the phase transition is sharp. She
showed that for a good choice of trees and car arrivals, a coupling
between the parking model and the Erdős--Rényi random graph model
enabled to study the critical window of the phase transition and
provided information about the geometry of the clusters of parked cars.
She established an unexpected link between the parking model and planar
maps by using a »last car« decomposition. This link has been opened
again in her contribution with Nicolas Curien (Contat and Curien, 2023).

In her initial follow-up work (Contat, 2022) she extended the results of
Curien and Hénard on general Bienaymé-Galton-Watson trees and allowed
different car arrival distributions depending on the vertex outdegrees.
She proved that this phase transition is sharp by establishing a large
deviations result for the flux of exiting cars.

In 2023, she has jointly with Nicolas Curien studied a combination of
parking on Cayley trees and a frozen modification of Erdős--Rényi random
graph model. Frozen here denotes slowing down the growth of components
which are not trees but contain cycles. They described phase transition
for the size of the components of parked cars using a modification of
the multiplicative coalescent which they called the frozen
multiplicative coalescent. They also studied geometry of critical parked
clusters. They relied on asymptotic results from Aldous (1997).~Derived
trees were very different from Bienaymé-Galton-Watson trees and should
converge towards the growth-fragmentation trees canonically associated
to the 3/2-stable process that already appeared in the study of random
planar maps. Already in her PhD work she pointed to some probable
connections to the study of random planar maps.

With David Aldous, Curien and Olivier Hénard (Aldous et al., 2023) she
studied parking on a infinite binary tree. Extensions of the parking
problem to binary and ordered trees have been already pointed to by
Lackner and Panholzer. Let \((A_{u}:u \in \mathbb{B)}\) be i.i.d.
non-negative integers that we interpret as car arrivals on the vertices
of the full binary tree \(\mathbb{B}\). It is known that the parking
process on \(\mathbb{B}\) exhibits a phase transition in the sense that
either a finite number of cars do not manage to park in expectation
(subcritical regime) or all vertices of the tree contain a car and
infinitely many cars do not manage to park (supercritical regime). They
characterized those regimes in terms of the law of \(A\) in an explicit
way and studied in detail the critical regime and the phase transition,
turning out to be discontinuous.

In another paper she studied parking on trees with a random given degree
sequence and the frozen configuration model (Contat, 2023), reminding on
her joint paper with Curien. She established and proved a natural
coupling between the frozen configuration model and the parking process
on a tree with prescribed degree sequence and prescribed car arrivals.
She established a phase transition for such process, as follows. She
firstly assumes a sequence of random degree sequences
\(I^{(n)} = (I_{1}^{(n)},\ldots,I_{n}^{(n)})\) and
\(A^{(n)} = (A_{1}^{(n)},\ldots,A_{n}^{(n)})\) such that for all \(n\),
the total in-degree \(\sum_{k = 1}^{n}I_{k}^{(n)}\) is equal to
\(n - 1\) and
\(\frac{1}{n}\sum_{k = 1}^{n}\delta_{(I_{k}^{(n)},A_{k}^{(n)})}\overset{\rightarrow}{n \rightarrow \infty}\lambda = \sum_{k \geq 0}^{}v_{k}\sum_{j \geq 0}^{}\mu_{(k),j}\delta_{(k,j)},\)
where the measure \(v = \sum_{k \geq 0}^{}{v_{k}\delta_{k}}\) is a
probability measure that we see as an offspring distribution and for all
\(k \geq 0\), the measure
\(\mu_{(k)} = \sum_{j \geq 0}^{}\mu_{(k),j}\delta_{j}\) is a probability
measure that represents the typical distribution of the car arrivals on
a vertex of out-degree \(k\). We assume that \(v\) has mean \(1\) and
finite non zero variance \(\Sigma^{2} \in (0,\infty)\), and for all
\(k \geq 0\), we assume that
\(\mu_{(k)} = \sum_{j \geq 0}^{}\mu_{(k),j}\delta_{j}\) has mean
\(m_{(k)} < \infty\) and finite variance \(\sigma_{(k)}^{2}\).

Theorem 2 (Contat, 2023): We assume the above assumptions with
\(\mathbb{E}_{\overline{\upsilon}}\lbrack m\rbrack \leq 1\) and
\(\mathbb{E}_{\upsilon}\lbrack m\rbrack \leq 1\). We also assume that
there esists a constant \(K\) such that \(m_{(k)} < K\) and
\(\sigma_{(k)}^{2} < K\) for all \(k \geq 0\). The parking process
undergoes a phase transition which depends on the sign of the quantity
\(\Theta = \left( 1 - \mathbb{E}_{\overline{\upsilon}}\lbrack m\rbrack \right)^{2} - \Sigma^{2}\mathbb{E}_{\upsilon}\lbrack\sigma^{2} + m^{2} - m\rbrack\).
More precisely, we have:

\[\frac{\varphi(T\left( I^{(n)} \right))}{n}{\overset{\mathbb{(P)}}{\rightarrow}}_{n \rightarrow \infty}C_{\lambda}\]

where \(C_{\lambda} = 0\) if and only if \(\Theta \geq 0\).

For the frozen configuration model, she proves the following local
scaling limit.

Proposition 5 (Contat, 2023): Suppose that
\(\mathbb{E}_{\overline{\upsilon}}\lbrack m\rbrack \leq 1\) and
\(\Theta \geq 0\). Under the above assumptions with
\(\mathbb{E}_{\upsilon}\lbrack m\rbrack \leq 1\), the frozen
configuration model converges Benjamini-Schramm quenched towards the
Bienaymé-Galton-Watson tree \(\mathcal{T}\) which is almost surely
finite.

In the above paper, Contat conjectures that no matter if one considers
strongly or weakly connected components, when the offspring distribution
of the tree or the car arrivals distributions have an infinite third
moments and a tail of order \(c \bullet k^{- \gamma}\) for some constant
\(c\) and some \(\gamma \in (3,4)\) when \(k\) goes to infinity, the
size of the components should be of order
\(n^{(\gamma - 2)/(\gamma - 1)}\) at criticality.

In a 2024 paper with Linxiao Chen (Chen and Contat, 2024), she studied
parking on supercritical geometric Bienaymé-Galton-Watson trees. She
provided a criterion to determine the phase of the parking process
(subcritical, critical, or supercritical) depending on the generating
function of \(\mu\). In a previous paper, Goldschmidt and Przykucki
(2019) proved that there are two possible regimes for the parking
process on the supercritical Bienaymé-Galton-Watson trees
\(\mathcal{T}\) depending on the two laws \(\mu\) and \(\upsilon\):

\begin{itemize}
\item
  Either \(\mathbb{E}\lbrack X\rbrack < \infty\) (subcritical regime)
\item
  Or \(X = \infty\) as soon as \(\mathcal{T}\) is infinite
  (supercritical regime)
\end{itemize}

In their paper, Chen and Contat mainly focus on the subcritical regime
of the above dichotomy. There main result in this case is the following.

Proposition 2 (\(F\)-characterization of the subcritical regime, Chen
and Contat, 2024): The law \(\mu\) is subcritical for the parking
process on a Bienaymé-Galton-Watson tree with geometric offspring
distribution with parameter \(q\) if and only if there exists a positive
solution \(p_{{^\circ}} > 0\) to the equation

\[\frac{1 - qp}{q} \bullet F\left( \frac{q(1 - q)}{(1 - qp)^{2}},1 \right) + p = 1.\]

As examples, they study geometric arrivals, Poisson arrivals and stable
cases, i.e. when the car arrivals distribution is non-generic.

In a 2025 paper, Contat and Lucile Laulin study parking on the random
recursive tree (Contat and Laulin, 2025). They prove that although the
random recursive tree has a non-degenerate Benjamini-Schramm limit, the
phase transition for the parking process appears at density \(0\). They
identify the critical window for appearance of a positive flux of cars
with high probability, which in the case of binary car arrivals happens
at density \({\log{(n)}}^{- 2 + o(1)}\) where \(n\) is the size of the
tree. Their work is the first that studied the parking process on trees
with possibly large degree vertices.

In her most recent joint paper with Curien (Contat and Curien, 2025),
she showed that critical parking trees conditioned to be fully parked
converge in the scaling limits towards the Brownian growth-fragmentation
tree, a self-similar Markov tree different from Aldous' Brownian tree
recently introduced and studied in Bertoin, Curien and Riera (2024).

Initially, she firstly derives the following result on asymptotics.

Corollary 3 (Universality of asymptotics for partition functions, Contat
and Curien, 2025). Under our standing assumptions on boundedness,
exchangeability, branching, aperiodicity for the flux, connectivity and
aperiodicity for the vertices, the functions
\(x \mapsto \left\lbrack y^{p} \right\rbrack F(x,y)\) for \(p \geq 0\)
have a common radius of convergence \(x_{cr} \in (0,\infty)\). For
\(x \in (0,x_{cr}\rbrack\), if \(y_{cr}^{x}\) is the radius of
convergence of \(y \mapsto F(x_{cr},y)\) then
\(y_{cr}^{x} \in (0,\infty)\) and \(F\left( x_{cr},y \right) < \infty\).
Furthermore, for each \(x \in (0,x_{cr}\rbrack\) there exists some
constants \(C^{x} > 0\) such that

\begin{itemize}
\item
  for \(x < x_{cr}\),
  \(W_{p}^{x} = \left\lbrack y^{p} \right\rbrack F(x,y)\sim C^{x} \bullet \left( y_{cr}^{x} \right)^{- p} \bullet p^{- \frac{3}{2}},\ as\ p \rightarrow \infty,\)
\item
  for
  \(x = x_{cr},\ W_{p}^{x_{cr}} = \left\lbrack y^{p} \right\rbrack F\left( x_{cr},y \right)\sim C^{x_{cr}} \bullet \left( y_{cr}^{x_{cr}} \right)^{- p} \bullet p^{- \frac{5}{2}},\ as\ p \rightarrow \infty.\)
\end{itemize}

They then study parking for self-similar Markov trees, defined as random
real rooted trees \(\mathcal{(T,}\rho)\) given with a decoration
\(g:\mathcal{T} \rightarrow \mathbb{R}_{+}\) starting from \(1\) at the
root \(\rho\) and positive on its skeleton. A special family of such
random trees
\(\mathfrak{T}_{\gamma} = (\mathcal{T}_{\gamma},g_{\gamma},\mu_{\mathcal{T}_{\gamma}})\)
was exhibited in relation with spectrally negative stable processes of
index \(\gamma \in (0,2)\). In particular,
\(\mathfrak{T}_{\frac{1}{2}}\) is nothing but a decorated version of the
famous Brownian continuous random tree of Aldous (Aldous, 1997), and
\(\mathfrak{T}_{\frac{3}{2}}\) is the Brownian growth-fragmentation tree
which already appeared inside the Brownian sphere/disk and was
conjectured to be the scaling limit of parking trees in Contat and
Curien (2023). They derive the following limit result.

Theorem 4 (Universal self-similar limits for the fully parked trees,
Contat and Curien, 2025). Under the standing assumptions on boundedness,
exchangeability, branching, aperiodicity for the flux, connectivity and
aperiodicity for the vertices, we have:

\begin{itemize}
\item
  When \(x < x_{cr}\) there exists some constant \(s^{x},v^{x} > 0\)
  such that
\end{itemize}

\[\left( s^{x} \bullet \frac{t}{p^{\frac{1}{2}}},\frac{\phi}{p},v^{x} \bullet \frac{\mu_{t}}{p} \right)under\ \mathbb{P}_{p}^{x}{\overset{(d)}{\rightarrow}}_{p \rightarrow \infty}\mathfrak{T}_{\frac{1}{2}},\]

\begin{itemize}
\item
  When \(x = x_{cr}\)
\end{itemize}

\[\left( s^{x_{cr}} \bullet \frac{t}{p^{\frac{3}{2}}},\frac{\phi}{p},v^{x_{cr}} \bullet \frac{\mu_{t}}{p^{2}} \right)under\ \mathbb{P}_{p}^{x}{\overset{(d)}{\rightarrow}}_{p \rightarrow \infty}\mathfrak{T}_{\frac{3}{2}},\]

the above convergence holds for the Gromov-Hausdorff-Prokhorov hypograph
convergence developed in Bertoin, Curien and Riera (2024).

\begin{enumerate}
\def\labelenumi{\arabic{enumi}.}
\setcounter{enumi}{4}
\item
  Future directions: metric topologies, connections to random planar
  maps and open problems of Lackner and Panholzer
\end{enumerate}

Study of parking problems offers a nice meeting bridge between
probability and combinatorics, including graph theory and discrete
mathematics. This can in future feature extensions in analysis, for
example in functional analysis by study in different metric topologies,
as well as complex analysis and geometry by extensions to the study of
planar maps with possible extensions to harmonic analysis,
representation theory and algebra.

At present, issues of studying parking problems in different metric
topologies remains largely unaddressed. Contat and Curien (2023) have
themselves pointed to appropriate contributions to follow in this sense:
Bhamidi, van der Hofstad and Sen (2018), Conchon-Kerjan and Goldschmidt
(2023) and Broutin, Duquesne and Wang (2018). This has been noted in the
context of extensions to the study of parking problems for car arrivals
with heavy tails, for example of a power-law distribution. It seems as
expected that scaling limit results, in particular when combined with
frozen modification of the Erdős-Rényi process would combine additive
and multiplicative coalescent, due to the connection with Gromov-weak
and Gromov-Hausdorff-Prokhorov topologies. Limiting results could relate
to inhomogenous continuum random trees as described in Bhamidi, van der
Hofstad and Sen (2018), and in this way resemble the results and
conjectures of Contat and Curien (2023). This would again relate study
of parking on random trees with the literature on random planar maps.

While Contat has in particular addressed the second and third open
problem noted by Lackner and Panholzer, her articles are far from
conclusive. In particular it would be very interesting to combine
additional extensions to random tree options and study in metric and
weak topologies of Gromov-Hausdorff-Prokhorov type. Dimensionality
issues have remained unaddressed in this line of research and one would
be tempted to ask if dimensions 3 and 4 would be special also for this
probabilistic problem, similar as they have proven to be in several
other cases in probability theory (Hutchcroft, 2025).

Problems 4, 5 and 6, noted by Lackner and Panholzer also remain
underaddressed and in need of further study in probability theory.
Quantities for tree and mapping parking functions such as total
displacement or individual displacements seem to lead to interesting
questions. Here, connections to results from queueing theory might merit
some interest and future possibilities -- these two areas seem to
feature a lot of possible resemblances. Study of defective parking
functions could be interesting to study in several contexts developed in
articles of Contat and extensions thereof. Enumeration problems for
restricted parking functions for trees or mappings could also be
interesting to study in present contexts developed and their above noted
extensions.

Additional possibilities for research could also be found in extensions
of the frozen modification of the Erdős-Rényi process, or even other
random graph possibilities, such as preferential attachment models.
Connections to present line of research in network archaeology would be
interesting to explore. Finally, connections to other present
probability and stochastic process research strands, in particular
interacting particle systems or even random matrix theory provide an at
present blank field of research.

References

\begin{enumerate}
\def\labelenumi{\arabic{enumi})}
\item
  Aldous, D. 1997. Brownian excursions, critical random graphs and the
  multiplicative coalescent. Ann. Probab. 25(2): 812-854 (April 1997).
  DOI: 10.1214/aop/1024404421.
\item
  Aldous, D., Contat, A., Curien, N., Hénard, O. 2023. Parking on the
  infinite binary tree. Probab. Theory Relat. Fields 187:481-504, DOI:
  10.1007/s00440-023-01189-6.
\item
  Armstrong, D., Reiner, V., Rhoades, B. 2015. Parking spaces. Adv.
  Math., 269:647706.
\item
  Berget, A., Rhoades, B. 2014. Extending the parking space. J. Combin.
  Theory Ser. A, 123:43-56.
\item
  Bertoin, J., Curien, N., Riera, A., 2024. Self-similar markov trees
  and scaling limits. arXiv preprint arXiv:2407.07888.
\item
  Bhamidi, S., van der Hofstad, R., Sen, S. 2018. The multiplicative
  coalescent, inhomogeneous continuum random trees, and new universality
  classes for critical random graphs. Probability Theory and Related
  Fields, 170:387-474.
\item
  Broutin, N., Duquesne, T., Wang, M. 2018. Limits of multiplicative
  inhomogeneous random graphs and L\(é\)vy trees, arXiv preprint
  arXiv:1804.05871.
\item
  Chen, L., Contat, A. 2024. Parking on supercritical geometric
  Bienaymé-Galton-Watson trees. arXiv preprint arXiv:2402.05612.
\item
  Conchon-Kerjan, G., Goldschmidt, C. 2023. The stable graph: the metric
  space scaling limit of a critical random graph with iid power-law
  degrees. Annals of Probability, 51(1):1-69. DOI: 10.1214/22-AOP1587.
\item
  Contat, A. 2022. Sharpness of the phase transition for parking on
  random trees. Random Structures and Algorithms, Volume 61, Issue 1,
  August 2022:84-100.
\item
  Contat, A. 2023. Parking on trees with a (random) given degree
  sequence and the Frozen configuration model. arXiv preprint
  arXiv:2312.04472.
\item
  Contat, A., Curien, N. 2023. Parking on Cayley trees and Frozen
  Erdős-Renyi. Annals of Probability, 51(6):1993-2055. DOI:
  10.1214/23-AOP1632.
\item
  Contat, A., Curien, N. 2025. Universality for catalytic equations and
  fully parked trees. arXiv preprint arXiv:2503.17348v1.
\item
  Contat, A., Laulin, L. 2025. Parking on the random recursive tree.
  arXiv preprint arXiv:2501.03195.
\item
  Cori, R., Le Borgne, Y. 2003. The sand-pile model and Tutte
  polynomials. Adv. in Appl. Math., 30(1-2):4452.
\item
  Foata, D., Riordan, J. 1974. Mappings of acyclic and parking
  functions. Aequationes Math., 10:10-22.
\item
  Goldschmidt, C., Przykucki, M. 2019. Parking on a random tree.
  Combinatorics, Probability and Computing, 28(1):23-45.
\item
  Haiman, M.D. 1994. Conjectures on the quotient ring by diagonal
  invariants. J. Algebraic Combin., 3(1):1776.
\item
  Hutchcroft, T. 2025. Critical cluster volumes in hierarchical
  percolation. Proc.Lond.Math.Soc. 130 (2025) 1, e70023.
\item
  Konheim, A.G., Weiss, B. 1966. An occupancy discipline and
  applications. SIAM J. Applied Math. 14:1266-1274.
\item
  Konvalinka, M., Tewari, V. 2021. Some natural extensions of the
  parking space. Journal of Combinatorial Theory, Series A, Volume 180,
  2021, 105394, ISSN 0097-3165, DOI: 10.1016/j.jcta.2020.105394.
\item
  Kung, J.P.S., Yan, C. 2003. Gončarov polynomials and parking
  functions. Journal of Combinatorial Theory Ser. A, 102:16-37.
\item
  Lackner, M.-L., Panholzer, A. 2016. Parking functions for mappings.
  Journal of Combinatorial Theory, Series A, 142:1-28.
\item
  Macdonald, I.G. 1995. Symmetric functions and Hall polynomials. Oxford
  Mathematical Monographs. The Clarendon Press, Oxford University Press,
  New York, second edition, 1995.
\item
  Pak, I.M., Postnikov, A.E. 1994. Resolvents for Sn-modules that
  correspond to skew hooks, and combinatorial applications. Funktsional.
  Anal. i Prilozhen., 28(2):7275.
\item
  Peterson, W.W. 1957. Addressing for random access storage. IBM J. Res.
  Develop., 1:130-146.
\item
  Stanley, R.P., Pitman, J. 2002. A polytope related to empirical
  distributions, plane trees, parking functions, and the associahedron.
  Discrete Comput. Geom., 27(4):603634.
\end{enumerate}

\end{document}